\documentclass{article}

\usepackage[para]{footmisc}

\DefineFNsymbols*{lamportnostar}[math]{\dagger\ddagger\S\P\|{\dagger\dagger}{\ddagger\ddagger}}
\setfnsymbol{lamportnostar}

\usepackage{authblk}
\usepackage{graphicx} 
\usepackage{sectsty}
\usepackage{etoolbox}
\usepackage{cite}
\usepackage{amsthm}
\usepackage{amsmath}
\usepackage{titlesec}
\usepackage{bm}
\usepackage [english]{babel}
\usepackage [autostyle, english = american]{csquotes}
\MakeOuterQuote{"}

\makeatletter
\def\CTEX@section@format{\Large\bfseries}
\makeatother

\patchcmd{\section}{\centering}{}{}{}

\theoremstyle{definition}

\theoremstyle{definition}
\newtheorem{lemma}{Lemma}[section]

\theoremstyle{definition}
\newtheorem{theorem}{Theorem}[section]

\title{
\textbf{Existence and Construction of a Gr{\"o}bner Basis for a Polynomial Ideal\footnote{Abstract of a paper presented at the workshop, \emph{Combinatorial Algorithms in Algebraic Structures}, Europäische Akademie, Otzenhausen, West Germany, Sept. 30 - Oct. 4, 1985.}}
}

\author{\emph{Deepak Kapur\\Paliath Narendran}}
\affil{\emph{Computer Science Research Branch}\\\emph{Corporate Research and Development}\\\emph{General Electric Company}\\\emph{Schenectady, NY 12345}}
\date{}

\begin{document}
\maketitle{}

\section{Introduction}
\noindent We show the existence of a Gröbner basis for a polynomial ideal over $R[x_1,\dots,x_n]$ 
assuming the existence of a Gröbner basis for every ideal over the underlying ring $R$. We 
also give an algorithm for computing a Gröbner basis assuming there is an algorithm to
compute a Gröbner basis for every ideal over $R$. We show how most of the properties that
hold for Gröbner bases for ideals over $R$ extend to Gröbner bases for polynomial ideas over 
$R[x_1,\dots,x_n]$. These results are structurally similar to Hilbert's Basis Theorem.
\\\\
A Gröbner basis of a polynomial ideal I can be defined in two ways: either as (a) any basis $B$
of $I$ such that every polynomial $p$ in $I$ reduces to 0 with respect to $B$, or as (b) any basis $B$ of $I$
such that every polynomial in the ring has a unique normal form with respect to $B$. We call
bases satisfying (a) \emph{weak Gröbner bases} and those satisfying (b) \emph{strong Gröbner bases}. When
Buchberger introduced the concept of Gröbner basis [1] for the case when $R$ is a field, he defined them to satisfy (a).
Later [2] he showed that they also satisfy condition (b). 
Weak Gröbner bases for various kinds of rings have since been studied in [6, 7, 8, 9, 13].
In [10] we defined weak Gröbner bases over first-order rings and showed their importance in theorem proving in first-order predicate calculus. 
Algorithms for computing strong Gröbner bases have been given in [4] (for reduction rings) and [5] (for Euclidean rings). 
In this paper we unify these two streams and show how to construct weak (strong) Gröbner bases over $R[x_1,\dots,x_n]$ provided weak (strong) Gröbner bases can be constructed over the underlying structure $R$.

\section{Assumptions about the Underlying Ring R}
\noindent Let $R$ be a commutative ring with 1.
We assume that for every ideal over $R$, there exists a Gröbner basis for the ideal. 
Let $\to$ stand for a \emph{reduction} relation induced by a basis on $R$.
\small
We will assume the following properties of $\to$:
\normalsize

\begin{itemize}
    \item[(i)] A non-zero element $a$ in $R$ \emph{reduces to 0} by another element $b$ if $a$ is a multiple of $b$.
    \item[(ii)] If $a\to c$, then there must exist $b_1,\dots,b_k\in R,k\geq1$, such that $a-c\in(b_1,\dots,b_k);$ we then say that \emph{a reduces to c using} $\{b_1,\dots,b_k\}$ in $R$. Usually $k=1$; in that case, we say that $a\to c$ by $b_1$.
    \item[(iii)] If an element $a$ reduces using $\{b_1,\dots,b_k\}$ and each $b_i$ reduces by $\{c_1,\dots,c_j\}$, then $a$ also reduces by $\{c_1,\dots,c_j\}$. In particular, if $a$ reduces by $b$ and $b$ reduces by $c$, then $a$ also reduces by $c$.
\end{itemize}

\begin{definition}
    A basis $B$ of an ideal $I$ over $R$ is a \emph{weak Gröbner basis} iff every element in the ideal has 0 as a normal form using $B$.
\end{definition}
\begin{definition}
    A basis $B$ of an ideal $I$ over $R$ is a \emph{strong Gröbner basis} iff every element in $R$ has a unique normal form using $B$.
\end{definition}

\noindent A well-founded reduction relation $\to$ induces a well-founded ordering on $R$ as follows: 
$a<c$ iff there exists a reduction sequence $c\to c_1\to\dots\to a$. Elements of $R$ in normal form wrt $\to$ are the minimal elements wrt $<$.

\noindent We will assume that we know whether the definition of reduction of elements of $R$ wrt other elements in $R$ admits weak Gröbner bases or strong Gröbner bases.
Let Gröbner be a function which takes a finite basis $B$ of an ideal over $R$ as the input and gives a Gröbner basis of the ideal as output; obviously (Gröbner(B))=(B).

\section{Preliminaries}
\noindent Let $R[x_1,\dots,x_n]$ be a polynomial ring over $R$.
We assume a total ordering on indeterminates which can be extended in many different ways to a total ordering on terms [4,5,12].
The ordering on $R$ induced by $\to$ and the ordering on terms can be used to define an ordering on monomials and polynomials. 
Given a polynomial $p$, let $HD(p), HT(p),\text{ and }HC(p)$ stand for its head-monomial, head-term, and head-coefficient, respectively. 
Let $Rest(p)$ be $p-HD(p)$. We assume that the functions $HD,HT,\text{ and }HC$ can be used on a set of polynomials also; for instance, $HD(I)=\{m\;|\; m\text{ is a head-monomial of a polynomial in }I\}$.
\clearpage 
\noindent Let $B$ be a finite set of polynomials in $R[x_1,\dots,x_n];B=\{p_1,\dots,p_m\}$. Let $I=(p_1,\dots,p_m)$.
\\\\
Using a well-founded relation $\to$ on $R$, we can define a well-founded reduction relation on polynomials in $R[x_1,\dots,x_n]$; we will denote it also by $\to$.
A polynomial $p\to q$ using another polynomial $p_1$ iff:
\begin{itemize}
    \item[(i)] There is a monomial $ct$ in $p$ i.e., $p=ct+p'$, such that $HT(p_1)$ divides $t$, i.e., $t=t'\cdot HT(p_1)$, and 
    \item[(ii)] There exists a $d$ such that $c\to d$ using $HC(p_1)$ and $q=p-k\cdot t'\cdot p_1$, where $c=k\cdot HC(p_1)+d$. 
\end{itemize}
The above definition extends to the case when more than one polynomial is used for reduction. 
Weak and strong Gröbner bases can be defined for polynomial ideals in the same way as definitions given above. 
\\\\
Let $MinMon(M)=\{m\;|\;m$ is a monomial in $M$ and there is no other monomial $m'$ in $M$ that divides $m\}.$
Similarly, let $MinTerm(T)=\{t\;|\;t$ is a term in $T$ and there is no $t'$ in $T$ that divides $t\}.$
It follows easily from Hilbert's Basis Theorem as well as from Dickson's Lemma that for any set of terms $T$, $MinTerm(T)$ is finite.

\section{Properties of Monomials in \boldmath$R[x_1,\dots,x_n]$}
\noindent Given a term $t$, let $C(t,I)=\{c\;|\; ct\text{ in }HD(I)\}$. 
It is easy to see that $C(t,I)$ is an ideal over $R$; further, for any two distinct terms $t_1$ and $t_2$, such that $t_1$ divides $t_2$, $C(t_1,I)\subseteq C(t_2,I)$. For every term $t$ appearing in some monomial in $MinMon(HD(I))$, i.e., $t$ in $HT(MinMon(HD(I)))$, let $D(t,I)$ be the set of all coefficients of $t$ in $MinMon(HD(I))$.
Let $Divisors(t)$ be the set of all terms including $t$ which can divide $t$. Note that for some $t$, $C(t,I)$ may not be empty, whereas $D(t,I)$ may be empty.

\begin{lemma}
    For each $c$ in $C(t,I)$, there exists $t'$ in $Divisors(t)$ and an element $d$ in $D(t',I)$ \small such that $d$ divides $c$.
\end{lemma}
\normalsize
\begin{lemma}
    The set $HT(MinMon(HD(I)))$ is finite.
\end{lemma}
\begin{lemma}
    For any term $t$, $C(t,I)$ is the same as the ideal generated by the union of $D(t',I)$ for all $t'$ in $Divisors(t)$.
\end{lemma}
\begin{lemma}
    $G=\bigcup\left(\text{Gröbner}(D(t',I)),t'\in Divisors(t)\right)$ is a Gröbner basis for $C(t,I)$.
\end{lemma}

\clearpage

\section{Existence of Gröbner basis for an Ideal over \boldmath$R[x_1,\dots,x_n]$}
\noindent Given an ideal $I$, we construct its Gröbner basis $GB$ as follows: 
For each $t$ in $HT(MinMon(HD(I)))$, let $G(t,I)$ be a Gröbner basis of the ideal generated by $D(t,I)$ over $R$; i.e., $G(t,I)=\text{Gröbner}(D(t,I))$;
for each element $g$ in $G(t,I)$, include in $GB$ a minimal polynomial in $I$ with $g\cdot t$ as its head-monomial.

\begin{theorem}
    The set $GB$ is finite.
\end{theorem}
\begin{theorem}
    For every polynomial $p$ in $I$, $p$ reduces to 0 using $GB$.
\end{theorem}
\begin{theorem}
    If $R$ admits strong Gröbner bases, then $GB$ is also a strong Gröbner basis.
\end{theorem}
\noindent Is such a $GB$ unique for an ideal $I$? 
No, because Gröbner bases for the same ideal may differ modulo units.
However, if we assume that $R$ admits strong Gröbner bases and such a Gröbner basis is unique for every ideal, then $GB$ is unique for a particular term ordering.
\begin{theorem}
    If $R$ admits stronger Gröbner bases and every ideal has a unique Gröbner basis, then every ideal over $R[x_1,\dots,x_n]$ also has a unique Gröbner basis.
\end{theorem}

\section{An Algorithm for Computing a Gröbner Basis}
\noindent We assume that (i) there is a way to compute a Gröbner basis $G$ from a basis $(c_1,\dots,c_k)$ over $R$;
furthermore, it is also possible to compute representations of $c_1,\dots,c_k$ in terms of the elements of $G$ as well as representations of every element $g$ in $G$ in terms of $c_1,\dots,c_k$;
these constructions are indeed possible from the construction of a Gröbner basis as illustrated by Buchberger in [12]. In addition, (ii) it should also be possible to generate a finite basis for a module of solutions to linear homogenous equations over $R$.
\\\\
{\bf 6.1 Critical Pairs}
\\

\noindent
For any finite subset $F$ of $B$, we define, 
two types of \emph{critical pairs} among polynomials in the basis. called \emph{G-polynomials} and \emph{M-polynomials}, (to stand for Gröbner-polynomials and Module-polynomials, respectively) as follows:
Let F consist of the polynomials
\[p_1=c_1\cdot t_1+r_1\;,\;p_2=c_2\cdot t_2+r_2,\;\dots\;,\;p_j=c_j\cdot t_j+r_j\]
where $HD(p_i)=c_i\cdot t_i,\;HT(p_i)=t_i,\;HC(p_i)=c_i,$ and $Rest(p_i)=r_i$.
\\Let $t=lcm(t_1,t_2,\dots,t_j)=s_1\cdot t_1=s_2\cdot t_2=\dots=s_j\cdot t_j$ for some terms $s_1,\dots,s_j$.
\clearpage
\noindent\textbf{G-Polynomial}: Let $\{g_1,\dots,g_m\}$ be a Gröbner basis of $\{c_1,\dots,c_j\}$; each $g_i$ can be written (represented) as $g_i=h_{i,1}\cdot c_1+h_{i,2}\cdot c_2+\dots+h_{i,j}\cdot c_j$ for some $h_{i,1},\dots,h_{i,j}$. Then, for each $i$, $q_i=h_{i,1}\cdot s_1\cdot p_1+h_{i,2}\cdot s_2\cdot p_2+\dots+h_{i,j}\cdot s_j\cdot p_j$ is a G-polynomial. Note that the head-monomial of $q_i$ is $g_i\cdot t$.
\\\\\textbf{M-Polynomials}: Consider the module $M=\{\langle a_1,\dots,a_j\rangle\;|\;a_1\cdot c_1+\dots+a_j\cdot c_j=0\}$.
Let the vectors $A_1,A_2,\dots,A_v$ represent its basis, where for every $i$, $A_i=\langle b_{i,1},\dots,b_{i,j}\rangle$ for some $b_{i,1},\dots,b_{i,j}$ in $R$. 
Then, for each $i$, \\$q_i=b_{i,1}\cdot s_1\cdot p_1+b_{i,2}\cdot s_2\cdot p_2+\dots + b_{i,j}\cdot s_j\cdot p_j$ is an M-polynomial.
Note that the head-term of $q_i$ is less than $t$ in the term ordering, since the monomials involving $t$ get cancelled in the summation. 
\\\\We say that \emph{q is a G-polynomial} (respectively, \emph{M-polynomial})\emph{ of B} if $q$ is a G-polynomial (respectively, M-polynomial) of some finite subset $F$ of $B$.

\begin{theorem}
    Let $B=\{p_1,\dots,p_m\}$ be a set of polynomials such that all the G-polynomials and M-polynomials of $B$ reduce to 0 using $B$. Then every polynomial in $(B)$ reduces to 0 using $B$, thus implying that $B$ is a weak Gröbner basis.
\end{theorem}
\begin{theorem}
    Let $B$ be a basis as stated in Theorem 6.1 above. If the function Gröbner gives a strong Gröbner basis on $R$, then $B$ is a strong Gröbner basis.
\end{theorem}
\noindent Furthermore, if for every ideal in $R$, the function Gröbner gives a unique strong Gröbner basis, then from Theorem 5.4 and Theorems 6.1 and 6.2, it follows that $B$ is also a unique strong Gröbner basis.
\\\\
The completion algorithm for obtaining Gröbner bases should be obvious by now: given a basis $B$, compute the G-polynomials and M-polynomials, and their normal forms. Add all non-zero normal forms to the basis and repeat this step until no new polynomails are added to the basis.
This process of adding new polynomials will terminate because of the finite ascending chain property of Noetherian rings.

\section{Instances of the Algorithm over Various Structures}
\noindent It can be shown that most of the known algorithms are instances of the above algorithm with additional properties assumed on the ring $R$.
For instance, when $R$ is a PID, then a Gröbner basis of a basis over $R$ can be obtained as the gcd of the elements in the basis. 
This can be computed by considering pairs of elements at a time. Further, a module basis of a finite set of non-zero elements in a PID can also be computed by considering pairs of elements.
For any two non-zero elements, $a$ and $b$, $\{\langle lcm(a,b)/a, -lcm(a,b)/b\rangle\}$ generates their module basis.
It can be shown that if G-polynomials and M-polynomials from every pair of polynomials in a basis $B$ reduce to 0, then G-polynomials and M-polynomials for every non-empty subset $B'$ of $B$ also reduce to 0.
So, it is sufficient to consider pairs of polynomials at a time for critical-pair computation. These considerations also apply to algorithms over fields and Euclidean rings since both are instances of PIDs.
\\\\
Zacharias [7], Trinks (as quoted in [8]), and Schaller [6] (as quoted in [4]) have given weak Gröbner basis algorithms for a ring $R$ on which the ideal membership problem is solvable as well as a finite set of generators for the solutions of linear equations in $R$ can be found algorithmically.
These algorithms can also be viewed as special cases of our algorithm where only M-polynomials are considered.

\section{References}
\vspace{0.25cm}
\noindent 1. Buchberger, B., \emph{An Algorithm for Finding a Basis for the Residue Class Ring of a Zero-Dimensional Ideal.} (in German) Ph.D Thesis, University of Innsbruck, Austria, 1965.
\vspace{0.1cm}
\\
2. Buchberger, B., "A Theoretical Basis for the Reduction of Polynomials to Canonical Forms," \emph{ACM SIGSAM Bulletin}, August, 1976, pp. 19-29.
\vspace{0.1cm}
\\
3. Bachmair, L., and Buchberger, B., "A Simplified Proof of the Characterization Theorem for Gröbner-Bases," \emph{ACM SIGSAM Bulletin}, Vol. 14, No. 4, 1980, pp. 29-34.
\vspace{0.1cm}
\\
4. Buchberger, B., "A Critical-Pair/Completion Algorithm in Reduction Rings," \emph{Proc. Logic and Machines: Decision Problems and Complexity} (eds. Borger, Hasenjaeger, Rodding), Springer Verlag LNCS 171, 1984, pp. 137-161.
\vspace{0.1cm}
\\
5. Kandri-Rody, A., and Kapur, D., "An Algorithm for Computing a Gröbner Basis of a Polynomial Ideal over a Euclidean Ring," \emph{TIS Report No. 84CRD045}, General Electric Corporate Research and Development, Schenectady, NY, Dec., 1984.
\vspace{0.1cm}
\\
6. Schaller, S., \emph{Algorithmic Aspects of Polynomial Residue Class Rings}, Ph.D Thesis, Computer Science Tech, Report 370, University of Wisconsin, Madison, 1979.
\vspace{0.1cm}
\\
7. Zacharias, G. \emph{Generalized Gröbner Bases in Commutative Polynomial Rings}, Bachelor Thesis, Lab. for Computer Science, MIT, 1978.
\vspace{0.1cm}
\\
8. Dubuque, W.G., Gianni, P., Trager, B., and Zacharias, G., "Primary Decomposition of Polynomial Ideal via Gröbner Bases," Extended Abstract, Unpublished Manuscript, 1985.
\vspace{0.1cm}
\\
9. Lauer, M., "Canonical Representatives for Residue Classes of a Polynomial Ideal," \emph{SYMSAC 1976}, pp. 339-345.
\vspace{0.1cm}
\\
\noindent 10. Kapur, D., and Narendran, P., "An Equational Approach to Theorem Proving in First-Order Predicate Calculus." \emph{7th International Conf. on Artificial Intelligence}, Los Angeles, Calif., August, 1985.
\vspace{0.1cm}
\\
11. van der Waerden, B.L., \emph{Modern Algebra}, Vol. I and II, Fredrick Ungar Publishing Co., New York, 1966.
\vspace{0.1cm}
\\
12. Buchberger, B., "Gröbner Bases: An Algorithmic Method in Polynomial Ideal Theory," to appear in \emph{Recent Results in Multidimensional Systems Theory} (ed. N.K. Bose), Reidel, 1985.
\vspace{0.1cm}
\\
13. Pan, L., "On the D-Bases of Ideals in Polynomial Rings over Principal Ideal Domains," \emph{This Workshop Proceedings}.
\clearpage 
\end{document}